\documentclass[10pt]{article}

\usepackage{amsmath}
\usepackage{amssymb}
\usepackage{indentfirst}
\usepackage{graphics} 
\usepackage{color}

\makeatletter

\@addtoreset{equation}{section}
\makeatother
\def\<{\langle}
\def\>{\rangle}

\newtheorem{lem}{Lemma}[section]
\newtheorem{theo}{Theorem}[section]
\newtheorem{rem}{Remark}[section]

\makeatletter
   
   \@addtoreset{equation}{section}
\makeatother

\begin{document}
\title{\bf Fast energy decay for wave equation with\\ a monotone potential and an effective damping}
\author{Xiaoyan Li\thanks{xiaoyanli@hust.edu.cn} \\ {\small School of Mathematics and Statistics}\\ {\small Huazhong University of Science and Technology} \\ {\small Wuhan, Hubei 430074, PR China}\\ and\\Ryo Ikehata\thanks{Corresponding author:ikehatar@hiroshima-u.ac.jp} \\ {\small Department of Mathematics, Division of Educational Sciences}\\ {\small Graduate School of Humanities and Social Sciences} \\ {\small Hiroshima University} \\ {\small Higashi-Hiroshima 739-8524, Japan}}
\date{}
\maketitle
\begin{abstract}
We consider the total energy decay of the Cauchy problem for wave equations with a potential and an effective damping. We treat it in the whole one-dimensional Euclidean space ${\bf R}$. 
Fast energy decay like $E(t)=O(t^{-2})$ is established with the help of potential. The proofs of main results rely on a multiplier method and modified techniques adopted in \cite{II}. 
\end{abstract}
\section{Introduction}
\footnote[0]{Keywords and Phrases: wave equation; one-dimensional space; potential; space-dependent damping; multiplier method; energy decay.}
\footnote[0]{2000 Mathematics Subject Classification. Primary 35L05; Secondary 35B40, 35B45.}

%
%
%
%
%
%


We consider the Cauchy problem for wave equation with a general potential and a space-dependent damping in the  one dimensional Euclidean space:
\begin{align} \label{c1}
\left \{
\begin{aligned}
&u_{tt}(t,x) - u_{xx}(t,x) + V(x)u(t,x) + a(x)u_{t}(t,x) 
= 0, ~~(t,x)\in (0,\infty)\times {\bf R},	\\
&u(0,x)= u_{0}(x), ~~x \in {\bf R},\\
&u_{t}(0,x)= u_{1}(x), ~~x \in {\bf R},
\end{aligned}	
	\right.
\end{align}
where the initial data $u_0(x)$ and $u_1(x)$ are chosen from $H^1(\bf R)$ and $L^2(\bf R)$ respectively. Both initial values and solutions shall all take real values. In addition, we assume that there exists $R>0$ such that
\begin{equation}\label{c23}
	\text{supp}\,u_0 \cup \text{supp}\, u_1 \subset  B_R:=\{x : |x| \leq R \}.	
\end{equation}

Throughout this paper, we define
\[u_{t}=\frac{\partial u}{\partial t},\quad u_{tt}=\frac{\partial^2 u}{\partial t^2},\quad 
u_{x}=\frac{\partial u}{\partial x},\quad
 u_{xx}=\frac{\partial^2 u}{\partial x^2}.\]
For convenience, we denote the usual  $L^q({\bf R})$-norm ($q = 2,\infty$) by 
$\| \cdot\|_q$. In particular, the $L^2({\bf R})$-norm is 
denoted by $\| \cdot\|$. The total energy $E(t)$ of the solution $u(t,x)$ to problem (1.1) is defined by
\begin{equation}
	E(t)=\frac{1}{2}(\| u_t(t,\cdot)\|^2+\| u_x(t,\cdot)\|^2+\|\sqrt{V(\cdot)}u(t,\cdot)\|^2 ).
\end{equation}
Furthermore, we define the inner product of $L^2({\bf R})$ by
\[
(u,v)=\int_{\bf R} u(x) v(x) dx.
\]
The energy decay of wave equation with local damping has been studied by many scholars. Let us  start with the initial-boundary value  problem 
\begin{align} \label{c34}
	\left \{
	\begin{aligned}
		&u_{tt}(t,x) - \Delta u(t,x) + a(x)u_{t}(t,x) 
		= 0, ~~(t,x)\in (0,\infty)\times {\Omega},	\\
		&u(0,x)= u_{0}(x), ~u_{t}(0,x)= u_{1}(x), ~~x \in {\Omega},\\
		&u(t,x)=0,~~x \in {\partial \Omega},
	\end{aligned}	
	\right.
\end{align}
where $\Omega = {\bf R}^{n}\setminus \bar{{\cal O}} \subset {\bf R}^{n}$ is a smooth exterior domain. For the case of  effective damping near infinity, that is, the damper $a(x)$ behaves  as
\begin{equation*}
a(x)\geq 0 ~~\text{in}~~\Omega, ~~~~~a(x)\geq \varepsilon_1>0~~  \text{for}~ |x|>\text{constant},
\end{equation*}
Nakao \cite{Nakao-1} and Ikehata \cite{I} provided the decay estimates of the total energy. Particularly, in \cite{I}, by special multiplier method developed in \cite{IM}, more faster energy decay such as $E(t)=O(t^{-2})$ was derived for \eqref{c34} under the star-shaped obstacle  ${\cal O}$ (${\cal \bar O}:={\bf R}^n  \setminus \Omega $) and weighted initial data conditions. Later on, some 
precise results  were obtained in \cite{A} under the so called GCC (Geometric Control Condition) assumption and without assuming star-shaped obstacle ${\cal O}$. The authors in \cite{A} showed that 
\[
E(t)=O(t^{-\min\{{\frac{3n}{4}, 1+\frac{n}{2}}\}})
\]
 for $n\geq 2$.  A generalization of \cite{I} was also done in \cite{D} by removing a geometrical condition assumed on the obstacles. Since mentioned results rely on the Poincare inequality and/or the Hardy inequality, only exterior domain and/or  high dimension cases can be treated. To overcome these obstacles, a novel method was developed in \cite{IL}, which employs potential term to compensate for the lack of the Poincar\'e and Hardy inequalities in the whole one dimensional Euclidean space. This idea has its origin in \cite{I-1}. For the whole high dimensional Euclidean space ${\bf R}^n (n \geq3)$, similar conclusions to \cite{I} about the total energy can be obtained. In addition, there are  many interesting results about effective damping near infinity. For example, in \cite{Zua}, the exponential decay of the total energy for Klein-Gordon type wave equation  was obtained by weighted energy method.
 The diffusion phenomena of wave equation with (asymptotically) periodic and/or constant damping was investigated in \cite{OZP}, \cite{J} and \cite{Nishihara}. For time-space dependent damper $a(t, x)$ including a constant damping, the decay properties of energy and $L^{2}$-norm of solution were studied in \cite{Mo}, \cite{M} and \cite{M-2}. Another generalization of \cite{I}  in noncompact Riemannian manifold was also considered by \cite{Z}.

For the case of degenerating damping near infinity, that is 
\[ a(x)\geq 0 ~~\text{in}~~\Omega, ~~~~~a(x)\to 0~~ \text{for}~~|x| \to \infty, 
\]
there are also many results for \eqref{c34} with $\Omega={\bf R}^n $. The most commonly chosen form of $a(x)$ is given by
\begin{equation}\label{c35}
\frac{a_1}{(1+|x|)^{\alpha}}\leq a(x)\leq \frac{a_2}{(1+|x|)^{\alpha}}, ~~~~(\alpha>0,~~x\in {\bf R}^n).	
\end{equation}
When $\alpha=1$, that is the so-called critical damping, the authors in \cite{ITY} showed the fact that if $1 < a_1 < n$, then
$E(t) = O(t^{-a_1})$; while if $a_1\geq n$, it holds that $E(t) = O(t^{-n+\delta})$ with  $\delta>0$ small enough. The authors in \cite{ITY} derived such results under the compact support condition on the initial data. It should be mentioned that very recently, an interesting result has been introduced by Sobajima \cite{S} in the case of $\alpha = 1$ (critical damping case), which completely remove the compactness of support of the initial data that assumed in \cite{ITY}, and furthermore achieve applications to non-linear problems of \eqref{c34}. In \cite{KNW}, the exponential decay of the total energy to problem \eqref{c34} was obtained for critical damping if the initial data were taken as some special form.
When $\alpha\in [0,1)$, \eqref{c35} is called as sub-critical damping. In this case, it was proved in \cite{TY} that the
energy of solutions to problem \eqref{c34} decays at a polynomial rate $t^{-(n-\alpha)/(2-\alpha)-1+\delta}$ for small enough $\delta>0$ (see also \cite{RTY, RTY-2} for higher order energy decay).  It should be mentioned that in the case of $n = 1$, we see $\frac{n-\alpha}{2-\alpha} + 1 < 2$. As a decay rate $2$ is a key number in our paper. In \cite{W-2}, a large time behavior of solutions and total energy to the wave equation with effective damping and absorbing nonlinearity is deeply studied in ${\bf R}^{n}$ with $n \geq 1$ for some weighted initial data. While, when $\alpha>1$, we call \eqref{c35} as super-critical damping. 
Mochizuki \cite{M-2} showed a non-decay in general property of the total energy for super-critical damping. In the case of super-critical damping with $\alpha > 1$, the total energy is generally non-decaying as is pointed out by \cite{M-2}, but it can be seen from \cite{BR} that the local energy does indeed decay with some rate. What's more, there are many results about so-called diffusion phenomena for an effective damping near infinity, that is, the solution to \eqref{c34} with $\Omega={\bf R}^n$ is approximated by a solution to the corresponding parabolic problem (\cite{SW-2}). We refer the interested readers to \cite{S-1, SW, SW-1, RTY-3, Nishi, W} for an additional topic of diffusion phenomenon. In this paper, "effective" damping is taken to mean the case where energy decay always occurs and therefore $\alpha \in [0,1]$.

Synthesis of the above researches, it seems that few results, such as $E(t)=O(t^{-2})$, are established for $n=1$.
It should be mentioned that the energy decay like $E(t)=O(t^{-2})$ for $n=1$ was obtained in \cite{IK} to problem \eqref{c34}. However, only half space case, that is $\partial \Omega = (0,\infty)$ in \eqref{c34}, is treated. So, whether one can obtain fast energy decay like $E(t)=O(t^{-2})$ for the whole one dimensional Euclid space seems like an open problem. In this paper, with the help of potential term $V(x)$, we derive the fast energy decay for effective damping $a(x)$ near infinity by a spacial multiply method. Basically one-dimensional Cauchy problems seem to be difficult because there are few useful tools compared with the higher dimensional case. In this paper, one breakthrough is presented.\\

In order to derive the main results, the following hypothesis are imposed for $a(x)$ and $V(x)$.

\noindent ({\bf A.1}) $a \in C({\bf R})$, and there exists two positive constants $a_1$ and $a_2$ such that
$$\frac{a_1}{(1+|x|)^{\alpha}}\leq a(x)\leq \frac{a_2}{(1+|x|)^{\alpha}}, ~~x\in {\bf R},$$
where $0 \leq \alpha \leq 1$. 
\\

\noindent ({\bf A.2}) $V \in C^{1}({\bf R})$ is a bounded function that satisfies 
$$V(x)>0,~~x V_{x}(x)\leq 0,~~~x\in{\bf R.}$$


With above preparations, our main results are stated as follows.\\

First we state fast energy decay property of the total energy appearing in the sub-critical damping case $\alpha < 1$. 
\begin{theo} \label{c21}
	Assume {\rm ({\bf A.1})} with $0 \leq \alpha<1$ and {\rm ({\bf A.2})}. If the initial data $[u_{0},u_{1}] \in H^{1}({\bf R})\times L^{2}({\bf R})$ satisfies \eqref{c23}, there exists a unique weak solution $u \in C([0,\infty); H^{1}({\bf R}))\cap C^{1}([0,\infty); L^{2}({\bf R}))$ to problem $\eqref{c1}$ satisfying $u(t,x) = 0$ for $\vert x\vert > R + t$ {\rm (}$t \geq 0${\rm )}, and 
	\[E(t)=O(t^{-2})\quad (t\to \infty).\]
\end{theo}

Next, we consider the critical damping case $\alpha = 1$. This case is rather complicated.

\begin{theo} \label{c21-1}
	Assume {\rm ({\bf A.1})} with $\alpha=1$ and {\rm ({\bf A.2})}. If the initial data $[u_{0},u_{1}] \in H^{1}({\bf R})\times L^{2}({\bf R})$ satisfies \eqref{c23}, there exists a unique weak solution $u \in C([0,\infty); H^{1}({\bf R}))\cap C^{1}([0,\infty); L^{2}({\bf R}))$ to problem $\eqref{c1}$ satisfying $u(t,x) = 0$ for $\vert x\vert > R + t$ {\rm (}$t \geq 0${\rm )}, and the following properties.\\	

\noindent{\rm (1)}~If~$0< a_1 \leq 2$, it holds that
\begin{equation}\label{c36}
E(t)=O(t^{-a_1 + \delta})\quad (t\to \infty)
\end{equation}
with small enough $\delta >0$.\\

\noindent{\rm (2)}~If~$a_1 >2$, it holds that
\begin{equation}\label{c37}
	E(t)=O(t^{-2})\quad (t\to \infty).
\end{equation}

\end{theo}
\noindent
{\bf Example.} Let $V_{0} > 0$. One can present $V(x) := V_{0}e^{-x^{2}}$, $V(x) := V_{0}(1+x^{2})^{-\frac{\mu}{2}}$ ($\mu > 0$) and $V(x) := V_{0}$ as examples. The last one corresponds to the so-called Klein-Gordon equation (cf. \cite{Zua}).
\begin{rem}{\rm As long as we properly check the unique existence of the weak solution, our result holds formally for the case $\alpha < 0$. We note that the coefficient $a(x)$ of the damping term is spatially unbounded in the negative $\alpha$ case (cf. \cite{IT}).}
\end{rem}
\begin{rem}{\rm In the case of $a(x) = a_{1} > 0$ and $V(x) = 0$, from \cite{M} one can know the total energy decay such that $E(t) = O(t^{-\frac{3}{2}})$ ($t \to \infty$). When we compare Theorem \ref{c21} with Matsumura's estimate, one can get faster decay rate such as  $E(t) = O(t^{-2})$ ($t \to \infty$). The influence of potential $V(x)$ is strongly effective in our theory even in the case of weakly effective potential (rapidly decay potential) such as $V(x) := V_{0}e^{-x^{2}}$. It is still open to get  faster decay $E(t) = O(t^{-2})$ ($t \to \infty$) in the case of $V(x) = 0$.}
\end{rem}
\begin{rem}{\rm A found number $2$ in Theorem \ref{c21-1} about $a_{1}$ seems to be a threshold which divides the decay property into two parts \eqref{c36} and \eqref{c37}. \eqref{c36} may express a wave like property of the solution, and \eqref{c37} may imply diffusive aspect of the solution. These two types of properties of the solution are closely related to that of \cite{RTY} which treated higher dimensional case. From the observation in \cite[Remark 1.2]{IK} the critical number $2$ on the coefficient $a_{1}$ seems reasonable.}
\end{rem}
\begin{rem}{\rm In previous studies \cite{G, lai}, similar model to (1.1) was studied from the view point of critical exponent of the power of the nonlinearity, however, the authors in \cite{G, lai} did not treat the one dimensional case. Although we are dealing with a linear problem, it could be a milestone when dealing with future one-dimensional non-linear problems.}
\end{rem}
\begin{rem}{\rm We derive our results assuming condition \eqref{c23}, but the consideration of the case without that condition \eqref{c23} is so far unknown.}
\end{rem}
The remainder of this paper is organized as follows. In Sections 2 and 3 we give proofs of Theorems \ref{c21} and \ref{c21-1} by modifying a method, which has its origin in \cite{II}.


\section{Proof of  main results}

In this section, we prove our main results. Since the unique existence of the weak solution $u(t,x)$ satisfying the finite speed of propagation property to problem (1.1) is a standard argument (cf. \cite{ikawa}) it suffices to get only the desired decay estimates in each Theorems \ref{c21} and \ref{c21-1}. The argument developed in \cite{II, IK} is useful again.\\



A next energy identity will play a crucial role in our argument.
\begin{lem}
For the solution $u(t,x)$ to Cauchy  problem \eqref{c1}, it holds that
\begin{align} 
&\frac{1}{2} \frac{d}{dt}G(t)+\frac{1}{2}\int_{\bf R} F_1(t,x) |u_t(t,x)| ^2dx+\frac{1}{2}\int_{\bf R} F_2(t,x) |u_x(t,x)| ^2dx \notag \\
&+\frac{1}{2}\int_{\bf R} F_3(t,x) |u(t,x)| ^2dx+\int_{\bf R} F_4(t,x) u_x(t,x)u_t(t,x)dx \notag\\
&+\frac{1}{2}\int_{\bf R} \frac{\partial}{\partial x}K(t,x) dx=0, \label{c4}
\end{align}
where
\begin{align*}
	~~~~~~~~~~~~~G(t)=&\int_{\bf R}f(t)\big(|u_t(t,x)|^2+|u_x(t,x)|^2+V(x)|u(t,x)|^2\big)+2g(t)u(t,x)u_t(t,x)\\
	&+\big(g(t)a(x)-g_t(t)\big)|u(t,x)|^2+2h(t,x)u_t(t,x)u_x(t,x)dx,\qquad\qquad\qquad\qquad\qquad\qquad\qquad\qquad\\
F_1(t,x)&=2f(t)a(x)-f_t(t)-2g(t)+h_x(t,x),   \\  
F_2(t,x)&=2g(x)-f_t(t)+h_x(t,x), \vspace{0.5ex} \\	
F_3(t,x)&=g_{tt}(t)-g_t(t)a(x)-V(x)f_t(t)+2V(x)g(t)-V_x(x)h(t,x)-V(x) h_x(t,x),\vspace{1ex}\\
F_4(t,x)&=h(t,x)a(x)-h_t(t, x),\\
K(t,x)&=-2f(t,x)u_t(t,x)u_x(t,x)-2g(t)u(t,x)u_x(t,x)-h(t,x)|u_t(t,x)|^2\\
&~\,~~-h(t,x)|u_x(t,x)|^2+V(x)h(t,x)|u(t,x)|^2.
\end{align*}
Here, $f(t)$, $g(t)$ and $h(t,x)$ are all smooth functions, which will be determined later on.
\end{lem}
{\it Proof.} The proof can be done for the smooth solution $u(t,x)$ by density. The argument developed in \cite{N, Nakao-1} is helpful.   
To make the proof more clear, we divide the following statement into four steps.	

\underline{{\it step 1.}} Multiplying the both sides of \eqref{c1} by $f(t)u_t(t,x)$ yields
\begin{align}
0=&fu_tu_{tt}-fu_tu_{xx}+fVuu_t+fa|u_t|^2 \notag\\
=&\frac{1}{2} f \frac{\partial}{\partial t}|u_t|^2-f\frac{\partial}{\partial x}(u_tu_x)+\frac{1}{2}\frac{\partial}{\partial t}(f|u_x|^2)-\frac{1}{2}f_t|u_x|^2+\frac{1}{2}fV\frac{\partial}{\partial t}|u|^2+fa|u_t|^2  \notag\\
=&\frac{1}{2} \frac{\partial}{\partial t}f |u_t|^2-\frac{1}{2}f_t|u_t|^2-\frac{\partial}{\partial x}(fu_tu_x)+\frac{1}{2}\frac{\partial}{\partial t}(f|u_x|^2)-\frac{1}{2}f_t|u_x|^2  \notag\\
&+\frac{1}{2}\frac{\partial}{\partial t}(fV|u|^2)-\frac{1}{2}Vf_t|u|^2+fa|u_t|^2  \notag\\
=&\frac{1}{2} \frac{\partial}{\partial t}(f|u_t|^2+f|u_x|^2+fV|u|^2)+(fa-\frac{1}{2}f_t)|u_t|^2  \notag\\
&-\frac{1}{2}f_t|u_x|^2-\frac{1}{2}f_tV|u|^2-\frac{\partial}{\partial x}(fu_tu_x) \label{c2}
\end{align}

\underline{{\it step 2.}} Multiplying the both sides of \eqref{c1} by $g(t)u(t,x)$, we obtain
\begin{align}
	0=&guu_{tt}-guu_{xx}+gV|u|^2+gauu_t \notag\\
	=&g \frac{\partial}{\partial t}(uu_t)-g|u_t|^2-g\frac{\partial}{\partial x}(uu_x)+g|u_x|^2+gV|u|^2+\frac{1}{2}ga\frac{\partial}{\partial t}|u|^2  \notag\\
	=& \frac{\partial}{\partial t}(guu_t)-\frac{1}{2}\frac{\partial}{\partial t}(g_t|u|^2)+\frac{1}{2}g_{tt}|u|^2-g|u_t|^2-\frac{\partial}{\partial x}(guu_x) \notag\\
	& +g|u_x|^2+gV|u|^2+\frac{1}{2}\frac{\partial}{\partial t}(ga|u|^2 )-\frac{1}{2}g_t a |u|^2\notag\\
	=&\frac{1}{2}\frac{\partial}{\partial t}(2guu_t-g_t\vert u\vert^{2}+g a |u|^2 )
	-g|u_t|^2+g|u_x|^2 \notag \\
	&+(\frac{1}{2}g_{tt}+gV-\frac{1}{2}g_t a)|u|^2-\frac{\partial}{\partial x}(guu_x).
\end{align}

\underline{{\it step 3.}} Multiplying the both sides of \eqref{c1} by $h(t,x)u_x(t,x)$, we have
\begin{align}
	0=&hu_xu_{tt}-hu_xu_{xx}+hVu_xu+hau_x u_t \notag\\
	=&h \frac{\partial}{\partial t}(u_xu_t)-hu_{xt} u_t-\frac{1}{2}h\frac{\partial}{\partial x}|u_x|^2+\frac{1}{2}hV \frac{\partial}{\partial x }|u|^2+hau_xu_t\notag\\
	=&  \frac{\partial}{\partial t}(hu_xu_t)-h_tu_xu_t-\frac{1}{2}\frac{\partial}{\partial x}(h|u_t|^2)+\frac{1}{2}h_x|u_t|^2-\frac{1}{2}
	\frac{\partial}{\partial x}(h|u_x|^2)\notag\\
	&+\frac{1}{2}h_x |u_x|^2+\frac{1}{2}
	\frac{\partial}{\partial x}(hV|u|^2)-\frac{1}{2}(h_xV+hV_x)|u|^2+ hau_xu_t\notag\\
	=&\frac{\partial}{\partial t}(hu_xu_t)+\frac{1}{2}h_x|u_t|^2+\frac{1}{2}h_x |u_x|^2-\frac{1}{2}(h_xV+hV_x)|u|^2\notag \\
	&+\frac{1}{2} \frac{\partial}{\partial x}(-h|u_t|^2-h|u_x|^2+hV|u|^2)
  +(ha-h_t)u_xu_t. \label{c3}
\end{align}

\underline{{\it step 4.}} Adding these identities from \eqref{c2} to \eqref{c3} all together, and integrating it over ${\bf R}$, we will arrive at the desired identity \eqref{c4}.
\hfill
$\Box$
\vspace{0.5cm}\\
Due to compactness of support of the initial data and the finite speed of propagation property of waves, the solution to problem \eqref{c1} vanishes for large $|x| \gg 1$. Therefore, we have 
\begin{align} \label{c5}
	\int_{\bf R} \frac{\partial}{\partial x}K(t,x) dx=K(t,+\infty)-K(t,-\infty)=0.
\end{align}
While, it follows from the Young inequality that
\begin{align}
\int_{\bf R} F_4(t,x) u_x(t,x)u_t(t,x)dx\geq -&\frac{k}{2}\int_{\bf R} |h(t,x)|a(x)|u_t(t,x)|^2dx-\frac{1}{2k}\int_{\bf R} |h(t,x)|a(x)|u_x(t,x)|^2dx \notag \\
-&\frac{1}{2}\int_{\bf R} |h_t(t,x)||u_t(t,x)|^2dx-\frac{1}{2}\int_{\bf R} |h_t(t,x)| |u_x(t,x)|^2dx.
\label{c6}
\end{align}
Substituting \eqref{c5} and  \eqref{c6} to  \eqref{c4}, we get 
\begin{align}  \label{c11}
\frac{d}{dt}G(t)+\int_{\bf R} K_1(t,x) |u_t(t,x)| ^2dx+\int_{\bf R} K_2(t,x) |u_x(t,x)| ^2dx 
	+\int_{\bf R} F_3(t,x) |u(t,x)| ^2dx \leq 0,
\end{align}
where
\[ K_1(t,x)=F_1(t,x)-k|h(t,x)|a(x)-|h_t(t,x)|,
\]
\[ K_2(t,x)=F_2(t,x)-\frac{1}{k}|h(t,x)|a(x)-|h_t(t,x)|.
\]

Next, we specify the expressions of $f(t)$, $g(t)$ and $h(t,x)$ as 
\[
f(t)=\varepsilon_1 (1+t)^2,~~~~g(t)=\varepsilon_2 (1+t),~~~~h(t,x)=\varepsilon_3 (1+t)x \phi(x)\]
where
\[\phi(x) = \left\{
\begin{array}{ll}
	\displaystyle{1},&
	\qquad \vert x\vert \leq 1, \\[0.2cm]
	\displaystyle{\frac{1}{|x|}},& \qquad \vert x\vert \geq 1,
\end{array} \right. 
\]
and $\varepsilon_1$, $\varepsilon_2$ and $\varepsilon_3$ are some positive constants, which are determined later on. Note that the function $\phi(x)$ is Lipschitz continuous on ${\bf R}$.

With above preparations, we provide the estimates for $K_1(t,x)$ and $K_2(t,x)$ by the next lemma. It should be noted that these three positive constants $\varepsilon_1$, $\varepsilon_2$ and $\varepsilon_3$ play important roles 
 in the proof of following lemmas.

\begin{lem} \label{c9}
Suppose $a_1>2$ for $\alpha=1$  and $a_1>0$ for $0 \leq \alpha<1$ in \eqref{c35}. If all parameters $\varepsilon_1$, $\varepsilon_2$ and $\varepsilon_3$ are well-chosen,  for $t>t_0\gg 1$, it holds that 
\[
{\rm (\romannumeral1)} ~  K_1(t,x)\geq 0, ~~~x\in{\bf R},\quad\quad\quad\quad  \quad\quad    
{\rm (\romannumeral2)} ~  K_2(t,x)\geq 0,~~~x\in{\bf R}.
\]
\end{lem}
{\it Proof.}
We first divide the integral region  into two parts  $|x|\leq 1$ and $|x|> 1$, then we check ${\rm (\romannumeral1)}$ and {\rm (\romannumeral2)}.\\

\noindent ${\rm (\romannumeral1)}$  For the case of $|x|\leq1$, we have
\begin{align}
K_1(t,x)
=&2f(t)a(x)-f_t(t)-2g(t)+h_x(t,x)-k|h(t,x)|a(x)-|h_t(t,x)|	\notag \\
=&2\varepsilon_1(1+t)^2a(x)-2\varepsilon_1(1+t)-2\varepsilon_2(1+t)+\varepsilon_3(1+t)-k\varepsilon_3(1+t)|x|a(x)-\varepsilon_3|x| \notag\\
\geq & C_{\alpha}\varepsilon_1a_1(1+t)^2-2\varepsilon_1(1+t)-2\varepsilon_2(1+t)+\varepsilon_3(1+t)- k\varepsilon_3(1+t)a(x)-\varepsilon_3
\notag\\
=&(1+t)^2\big\{
C_{\alpha}\varepsilon_1 a_1-\frac{2\varepsilon_1}{1+t}-\frac{2\varepsilon_2}{1+t}+\frac{\varepsilon_3}{1+t}-\frac{k\varepsilon_3\|a\|_\infty}{1+t}-\frac{\varepsilon_3}{(1+t)^2}
\big\}
\end{align}
with some $\alpha$-dependent constant $C_{\alpha} > 0$, where one has just used the assumption {\bf (A.1)}.

 For the case of $|x|>1$, assumption {\bf (A.1)} and the finite speed of propagation property of the solution lead to

\begin{align}
	K_1(t,x)
	=&2f(t)a(x)-f_t(t)-2g(t)+h_x(t,x)-k|h(t,x)|a(x)-|h_t(t,x)|	\notag \\
	=&2\varepsilon_1(1+t)^2a(x)-2\varepsilon_1(1+t)-2\varepsilon_2(1+t)-k\varepsilon_3(1+t)a(x)-\varepsilon_3\notag\\
	\geq & 2\varepsilon_1 (1+t)^2\frac{a_1}{(1+|x|)^{\alpha}}-2\varepsilon_1(1+t)-2\varepsilon_2(1+t)-k\varepsilon_3(1+t)\frac{a_2}{(1+|x|)^{\alpha}}-\varepsilon_3
	\notag\\
	=&\frac{(1+t)^2}{(1+|x|)^{\alpha}}\big\{
	2\varepsilon_1 a_1-\frac{2\varepsilon_1(1+|x|)^{\alpha}}{1+t}-\frac{2\varepsilon_2(1+|x|)^{\alpha}}{1+t}
	-\frac{k\varepsilon_3 a_2}{1+t}-\frac{\varepsilon_3(1+|x|)^{\alpha}}{(1+t)^2} \big\}  \notag \\
	\geq &\frac{(1+t)^2}{(1+|x|)^{\alpha}}\big\{
	2\varepsilon_1 a_1-\frac{2\varepsilon_1(1+R+t)^{\alpha}}{1+t}-\frac{2\varepsilon_2(1+R+t)^{\alpha}}{1+t}
	\notag \\&
	-\frac{k\varepsilon_3 a_2}{1+t}-\frac{\varepsilon_3(1+R+t)^{\alpha}}{(1+t)^2} \big\}  
\end{align}
for large $t \gg 1$.\\

\noindent ${\rm (\romannumeral 2)}$  For the case of $|x|\leq1$, $K_2(t,x)$ satisfies
\begin{align}
	K_2(t,x)
	=&2g(t)-f_t(t)+h_x(t,x)-\frac{1}{k}|h(t,x)|a(x)-|h_t(t,x)|	\notag \\
\geq &2\varepsilon_2(1+t)-2\varepsilon_1(1+t)+\varepsilon_3(1+t)-\frac{1}{k}\varepsilon_3(1+t)a(x)-\varepsilon_3 \notag\\
	=&(1+t)\big\{
	2\varepsilon_2 -2\varepsilon_1+\varepsilon_3
	-\frac{\varepsilon_3\|a\|_\infty}{k}-\frac{\varepsilon_3}{1+t} \big\}.
\end{align}

 For the case of $|x|>1$, we have 
\begin{align}
	K_2(t,x)
	=&2g(t)-f_t(t)+h_x(t,x)-\frac{1}{k}|h(t,x)|a(x)-|h_t(t,x)|	\notag \\
	=&2\varepsilon_2(1+t)-2\varepsilon_1(1+t)-\frac{1}{k}\varepsilon_3(1+t)a(x)-\varepsilon_3 \notag\\ 	
   \geq &(1+t)\big\{
	2\varepsilon_2 -2\varepsilon_1
	-\frac{\varepsilon_3\|a\|_\infty}{k}-\frac{\varepsilon_3}{1+t} \big\}.
\end{align}

To guarantee the positivity of $K_1(t,x)$ and $K_2(t,x)$, our next task is to choose reasonable positive constants $k,~\varepsilon_1,~\varepsilon_2,$ and $\varepsilon_3$, which depend on the value of $\alpha \in [0,1]$. In this case, it is important to notice the fact that
\[\lim_{t \to \infty}\frac{(1+R+t)^{\alpha}}{1+t} = 0\quad(\alpha < 1), \quad \lim_{t \to \infty}\frac{1+R+t}{1+t} = 1.\]

\underline{Case for $0 \leq \alpha < 1$.}\\
For large $t\geq t_0\gg 1$, the following conditions are needed:
\begin{gather}
C_{\alpha}\varepsilon_1 a_1>0, \label{c19}\\
2\varepsilon_2-2\varepsilon_1+\varepsilon_3-\frac{ \varepsilon_3 \|a\|_\infty}{k}
> 0, \\
2\varepsilon_2-2\varepsilon_1-\frac{ \varepsilon_3 \|a\|_\infty}{k}
> 0. \label{c20}
\end{gather}
In fact, it is sufficient to choose $\varepsilon_1$ and $\varepsilon_2$ and large $k > 0$ satisfying
\begin{align*} 
	\left \{
	\begin{aligned}
		&\varepsilon_2>\varepsilon_1>0,
		\\
		&k>\frac{\varepsilon_3 \|a\|_\infty}{2(\varepsilon_2-\varepsilon_1)}.
	\end{aligned}	
	\right.
\end{align*}
In this case $\varepsilon_3 > 0$ can be chosen arbitrarily.\\

\underline{Case for $\alpha =1$.}\\ 
Additional conditions are needed.\\
\begin{gather} 
C_{\alpha}\varepsilon_1 a_1>0, \\	
2a_1\varepsilon_1 -2\varepsilon_1-2\varepsilon_2 > 0, \label{c7}\\
2\varepsilon_2-2\varepsilon_1+\varepsilon_3-\frac{ \varepsilon_3 \|a\|_\infty}{k}
> 0, \\
2\varepsilon_2-2\varepsilon_1-\frac{ \varepsilon_3 \|a\|_\infty}{k}
> 0.
 \label{c8} 
\end{gather}
To guarantee (2.17) and \eqref{c8}, it is necessary to set 
\begin{align*} 
	\left \{
	\begin{aligned}
		&\varepsilon_2>\varepsilon_1>0,
		\\
		&k > \frac{\varepsilon_3 \|a\|_\infty}{2(\varepsilon_2-\varepsilon_1)}.
	\end{aligned}	
	\right.
\end{align*}
Under this situation, to realize \eqref{c7} such that
\begin{equation}\label{Ike-1}
	(a_1-1)\varepsilon_1> \varepsilon_2>\varepsilon_1,
\end{equation}
we must choose $a_1>2$ in the case of $\alpha=1$. In this case $\varepsilon_3 > 0$ can be also chosen arbitrarily.

\hfill
$\Box$

We look at the bound for $F_{3}(t,x)$ in the next lemma. 

\begin{lem}\label{c10}
	Suppose $a_1>2$ for $\alpha=1$  and $a_1>0$ for $0 \leq \alpha<1$. If  the parameters $\varepsilon_1$,  $\varepsilon_2$ and $\varepsilon_3$ are well-chosen,  for $t>t_0\gg 1$, it holds that 
	\[
-F_3(t,x)\leq C a(x),~~~x\in{\bf R},
	\]
where $C > 0$ is a generous constant.
\end{lem}	
\noindent {\it Proof.}
By assumption {\bf (A.2)}, we have for a.e. $x \in {\bf R}$
\begin{align}
-F_3(t,x)=&-g_{tt}(t)+g_t(t)a(x)+V(x)f_t(t)-2V(x)g(t)+V_x(x)h(t,x)+V(x)h_x(t,x), \notag\\
\leq&\varepsilon_2 a(x)+2\varepsilon_1V(x)(1+t)-2\varepsilon_2V(x)(1+t)+\varepsilon_3V_x(x)x(1+t)\phi(x)+\varepsilon_3 V(x)(1+t)  \notag\\
\leq &\varepsilon_2a(x)+V(x)(1+t)(2\varepsilon_1-2\varepsilon_2+\varepsilon_3).
\end{align}
Let us choose parameters such as $2\varepsilon_1-2\varepsilon_2+\varepsilon_3\leq 0$, that is, $\varepsilon_3\leq 2\varepsilon_2-2\varepsilon_1$ to get the desired estimate. Here, it should be noted that we must choose additionally $2\varepsilon_2>2\varepsilon_1$. By the way, in the case of $\alpha=1$, by considering \eqref{Ike-1} we must choose simultaneously
\[
(a_1-1)\varepsilon_{1} > \varepsilon_2 > \varepsilon_1.
\]
Consequently, $a_1>2$ is again required. 
\hfill 
$\Box$\\

\noindent \underline{{\it Proof of Lemmas \ref{c9} and \ref{c10} completed.}}\, By summarizing the proof of Lemmas \ref{c9} and \ref{c10}, let us choose unified parameters $\varepsilon_j > 0$ ($j = 1,2,3$) and $k \gg 1$ for each $0 \leq \alpha \leq 1$.\\


\noindent{(1)}~\underline{In the case $0 \leq \alpha<1$}, we select $\mu>0$ such that
\[\varepsilon_1=\mu,~~~\varepsilon_2=3\mu,~~~ \varepsilon_3=\frac{1}{2}\mu,~~~ k>\frac{\|a\|_\infty}{8}.
\]

\noindent{(2)}~\underline{In the case $\alpha=1$}, we select $\lambda=a_1-2>0$ and  $\mu>0$  such that
\[\varepsilon_1=\mu,~~~\varepsilon_2=(1+\frac{\lambda}{2})\mu,~~~ \varepsilon_3=\frac{1}{4}\lambda \mu,~~~ k>\frac{\|a\|_\infty}{4}.
\]
With above parameters, the statements of Lemmas \ref{c9} and \ref{c10} are true.
\hfill
$\Box$\\

Thanks to  Lemmas \ref{c9} and \ref{c10}, \eqref{c11} can be rewritten as 
\[
\frac{d}{dt}G(t)\leq C\int_{\bf R} a(x)|u(t,x)|^2dx\quad (t \geq t_{0} \gg 1).
\]
Integrating it over $[t_0, t]$ yields
\begin{equation}\label{c12}
 G(t)\leq G(t_0)+C\int_{t_0}^t \int_{\bf R} a(x)|u(s,x)|^2dx ds.	
\end{equation}
We can estimate $G(t_0)$ by the Cauchy-Schwarz  inequality
\begin{align} 
G(t_0)\leq & 2\varepsilon_1(1+t_0)^2 E(t_0)+2\varepsilon_2(1+t_0)\|u(t_0, \cdot)\|\|u_t(t_0, \cdot)\|+\varepsilon_2 \|a(x)\|_\infty(1+t_0)\|u(t_0, \cdot)\|^2 \notag \\
&+2\varepsilon_3(1+t_0)\|u_t(t_0, \cdot)\| \|u_x(t_0, \cdot)\| 
:=C_{t_0}>0.
\end{align} \label{c13}
It follows from \eqref{c12} and the definition of $G(t)$ that
\begin{align}
	&f(t)E(t)+g(t)(u(t,\cdot), u_t(t,\cdot))+(h(t, \cdot)u_t(t,\cdot), u_x(t,\cdot)) \notag \\
	\leq& \frac{C_{t_0}}{2}+\frac{1}{2}\int_{\bf R}(g_t(t)-g(t)a(x))|u(t,x)|^2dx+\frac{C}{2}\int_{t_0}^t \int_{\bf R} a(x)|u(s,x)|^2dx ds \notag \\
	\leq& C_{t_0}+\varepsilon_2\Vert u(t,\cdot)\Vert^{2} + C\int_{t_0}^t \int_{\bf R} a(x)|u(s,x)|^2dx ds. \label{c14}
\end{align} 

Next, we need a crucial lemma to obtain main results, which can be derived by Ikehata \cite{I-1}. A role of potential $V(x)$ is crucial. We only write down its statement without proof. 

\begin{lem}\label{c17}\, Under the assumptions of Theorems {\rm \ref{c21}} and {\rm \ref{c21-1}} on the initial data, for the corresponding weak solution $u(t,x)$ to problem {\rm (1.1)-(1.2)}, it holds that
	\[\Vert u(t,\cdot)\Vert^{2} + \int_{0}^{t}\int_{{\bf R}}a(x)\vert u(s,x)\vert^{2}dxds \leq  C\left(\Vert u_{0}\Vert^{2} + \int_{{\bf R}}\frac{\vert u_{1}(x)+a(x)u_{0}(x)\vert^{2}}{V(x)}dx\right) =: CJ_{0}^{2}, \quad ~t \geq 0,\]
	where $C > 0$ is a generous constant.
\end{lem}

\begin{lem}
	\label{c18}\, Under the assumptions of Theorems {\rm \ref{c21}} and {\rm \ref{c21-1}}, for $t\geq t_0\gg 1$ it holds that
\begin{equation*}
	f(t)E(t)+(h(t, \cdot)u_x(t, \cdot), u_t(t,\cdot))\geq \frac{1}{2}f(t)E(t).
\end{equation*}
\end{lem}

\noindent {\it Proof.}
By the Young inequality, we have 
\begin{align}
&f(t)E(t)+(h(t, \cdot)u_x(t, \cdot), u_t(t,\cdot)) \notag\\ 
\geq &\frac{1}{2} \int_{\bf R} f(t)(|u_t(t,x)|^2+|u_x(t,x)|^2+V(x)|u(t,x)|^2)dx-\frac{1}{2} \int_{\bf R} |h(t,x)| (|u_t(t,x)|^2+|u_x(t,x)|^2)dx \notag	\\
\geq &\frac{1}{2} \int_{\bf R}(f(t)- |h(t,x)|) (|u_t(t,x)|^2+|u_x(t,x)|^2)dx+\frac{1}{4} \int_{\bf R} f(t)V(x) |u(t,x)|^2dx.  \label{c16}
\end{align} 
By choosing $t\geq t_0\gg1 $ large enough such that $$\frac{\varepsilon_3}{1+t}\leq \frac{\varepsilon_1}{2}, $$
we see
\begin{align} 
f(t)- |h(t,x)|
	\geq  (1+t)^2(\varepsilon_1- \frac{\varepsilon_3}{1+t})\geq \frac{\varepsilon_1}{2} (1+t)^2=\frac{1}{2}f(t). \label{c15}
\end{align} 

Substituting \eqref{c15} into \eqref{c16}, we can obtain the desired result.

\hfill
$\Box$

Finally, let us prove main results with the tools prepared above.

\noindent
\underline{ {\it Proof of Theorem \ref{c21}}}  ~~ 
Under the assumption of Theorem \ref{c21}, we can easily check 
\[\int_{{\bf R}}\frac{\vert u_{1}(x)+a(x)u_{0}(x)\vert^{2}}{V(x)}dx < \infty.\]
It follows from \eqref{c14}, Lemmas \ref{c17} and \ref{c18} that
\begin{align}
	f(t)E(t)\leq& 2g(t)\|u(t,\cdot)\| \|u_t(t, \cdot)\|+ C + CJ_{0}^{2}\notag \\
	\leq &C g(t)\sqrt{E(t)}J_{0} + C + CJ_{0}^{2}\notag \\
	\leq &\frac{CJ_{0}g(t)}{\sqrt{f(t)}}\sqrt{f(t)E(t)}+C + CJ_{0}^{2}\notag \\
	\leq &\frac{C^2}{2}J_{0}^{2}\frac{g(t)^2}{f(t)}+\frac{1}{2}f(t)E(t)+C + CJ_{0}^{2}
\end{align}
with some generous constant $C > 0$. Therefore,
\[	f(t)E(t)\leq C^2 J_{0}^{2}\frac{g(t)^2}{f(t)}+ C + CJ_{0}^{2}, \quad t \geq t_0,
\]
that is 
\[
E(t)\leq C^2 J_{0}^{2}\frac{g(t)^2}{f(t)^2}+\frac{C}{f(t)} + \frac{CJ_{0}^{2}}{f(t)},
\quad t \geq t_0, \]
which completes the proof.
\hfill
$\Box$
\\

\noindent
\underline{ {\it Proof of the case (2) of Theorem \ref{c21-1}}}  ~~ By similar arguments to Theorem \ref{c21}, we will get \eqref{c37} of Theorem \ref{c21-1} by assuming $a_{1} > 2$.
\hfill
$\Box$
\\


\section{Proof of the case (1) of Theorem \ref{c21-1}}

In a new chapter, we now prove the case (1) of Theorem \ref{c21-1}. The proof of the case (1) of Theorem \ref{c21-1} is similar to the previous part basically. So, we sketch out the main points in the following statements. 

This part starts with new definitions of $f(t)$, $g(t)$ and $h(t,x)$:
\[
f(t)=\varepsilon_1 (1+t)^\theta,~~~~g(t)=\varepsilon_2 (1+t)^{\theta-1},~~~~h(t,x)=\varepsilon_3 (1+t)^{\theta-1}x \phi(x),\]
where $\theta>0$ will be fixed later on.

Under these new definitions, we need to check  Lemmas  \ref{c9} and \ref{c10} with $0<a_1\leq 2$ and $\alpha=1$ for large $t >0$.\\

\noindent
{\bf Check of Lemma \ref{c9}}:~If  $|x|\leq1$, $K_1(t,x)$ satisfies
\begin{align}
	K_1(t,x)
=&2f(t)a(x)-f_t(t)-2g(t)+h_x(t,x)-k|h(t,x)|a(x)-|h_t(t,x)|	\notag \\
	=&2\varepsilon_1 (1+t)^\theta a(x)-\varepsilon_1 \theta (1+t)^{\theta-1}-2\varepsilon_2(1+t)^{\theta-1}+\varepsilon_3(1+t)^{\theta-1}\notag \\
	&-k\varepsilon_3(1+t)^{\theta-1}|x|a(x)-\varepsilon_3(\theta-1) (1+t)^{\theta-2}|x| \notag\\
	\geq &(1+t)^\theta \big\{
	\varepsilon_1 a_1C_{\alpha}-\frac{\theta \varepsilon_1}{1+t}-\frac{2\varepsilon_2}{1+t}+\frac{\varepsilon_3}{1+t}-\frac{k\varepsilon_3\|a\|_\infty}{1+t}-\frac{\varepsilon_3{(\theta-1)}}{(1+t)^2}
	\big\}
\end{align}
with some constant $C_{\alpha} > 0$. So, one can realize the positivity of $K_{1}(t,x)$ in $\vert x\vert \leq 1$ by taking large $t > 0$ for any $\theta > 0$.
If $|x|\geq1$, the finite speed of propagation property of the solution yields
\begin{align}
	K_1(t,x)
	=&2f(t)a(x)-f_t(t)-2g(t)+h_x(t,x)-k|h(t,x)|a(x)-|h_t(t,x)|	\notag \\
	=&2\varepsilon_1(1+t)^\theta a(x)-\theta \varepsilon_1(1+t)^{\theta-1}-2\varepsilon_2(1+t)^{\theta-1}-k\varepsilon_3(1+t)^{\theta-1}a(x)-\varepsilon_3{(\theta-1)} (1+t)^{\theta-2}\notag\\
%
	\geq &\frac{(1+t)^\theta}{(1+|x|)^{\alpha}}\big\{
	2\varepsilon_1 a_1-\frac{\theta \varepsilon_1(1+|x|)^{\alpha}}{1+t}-\frac{2\varepsilon_2(1+|x|)^{\alpha}}{1+t}
	-\frac{k\varepsilon_3 a_2}{1+t}-\frac{\varepsilon_3({\theta-1})(1+|x|)^{\alpha}}{(1+t)^2} \big\}  \notag \\
	\geq &\frac{(1+t)^\theta}{(1+|x|)^{\alpha}}\big\{
	2\varepsilon_1 a_1-\frac{\theta \varepsilon_1(1+R+t)^{\alpha}}{1+t}-\frac{2\varepsilon_2(1+R+t)^{\alpha}}{1+t}
	\notag \\&
	-\frac{k\varepsilon_3 a_2}{1+t}-\frac{\varepsilon_3({\theta-1})(1+R+t)^{\alpha}}{(1+t)^2} \big\}  
\end{align}
with $\alpha = 1$.

For $K_2(t,x)$, in the case of $\vert x\vert \leq 1$ it holds that
\begin{align}
	K_2(t,x)
	=&2g(t)-f_t(t)+h_x(t,x)-\frac{1}{k}|h(t,x)|a(x)-|h_t(t,x)|	\notag \\
	\geq &2\varepsilon_2(1+t)^{\theta-1}-\theta\varepsilon_1(1+t)^{\theta-1}+ \varepsilon_{3}(1+t)^{\theta-1}-\frac{1}{k}\varepsilon_3(1+t)^{\theta-1}a(x)-\varepsilon_3(\theta-1) (1+t)^{\theta-2}\notag\\
	\geq &(1+t)^{\theta-1}\big\{
	2\varepsilon_2 -\theta \varepsilon_1+\varepsilon_{3} 
	-\frac{\varepsilon_3\|a\|_\infty}{k}-\frac{\varepsilon_3{(\theta-1)}}{1+t} \big\}.
\end{align}

For $K_2(t,x)$, in the case of $\vert x\vert \geq 1$ it holds that
\begin{align}
	K_2(t,x)
	=&2g(t)-f_t(t)+h_x(t,x)-\frac{1}{k}|h(t,x)|a(x)-|h_t(t,x)|	\notag \\
	\geq &2\varepsilon_2(1+t)^{\theta-1}-\theta\varepsilon_1(1+t)^{\theta-1} -\frac{1}{k}\varepsilon_3(1+t)^{\theta-1}a(x)-\varepsilon_3(\theta-1) (1+t)^{\theta-2}\notag\\
	\geq &(1+t)^{\theta-1}\big\{
	2\varepsilon_2 -\theta \varepsilon_1
	-\frac{\varepsilon_3\|a\|_\infty}{k}-\frac{\varepsilon_3{(\theta-1)}}{1+t} \big\}.
\end{align}

To guarantee the positivity of $K_1(t,x)$ and $K_2(t,x)$, for large $t\geq t_0\gg 1$, the following conditions must be imposed:
\begin{gather*}
	C_{\alpha}\varepsilon_1 a_1>0, \\
	2a_1\varepsilon_1 - \varepsilon_1\theta-2\varepsilon_2 > 0,\\
	2\varepsilon_{2} - \theta\varepsilon_{1} + \varepsilon_{3}(1-\frac{\Vert a\Vert_{\infty}}{k}) > 0,\\
    2\varepsilon_2-\theta\varepsilon_1-\frac{ \varepsilon_3 \|a\|_\infty}{k}
	> 0, 
\end{gather*}
that is,
\begin{equation}\label{c29}
	k\geq \frac{\varepsilon_3 \|a\|_\infty}{2\varepsilon_2-\theta \varepsilon_1},~~ k > \Vert a\Vert_{\infty},~~ \frac{2a_1-\theta}{2}\varepsilon_1> \varepsilon_2>\frac{\theta}{2} \varepsilon_1,
\end{equation}
which implies  
\begin{equation}\label{c30}
	\frac{2a_1-\theta}{2}>\frac{\theta}{2}~~ \Rightarrow~~ \theta<a_1.
\end{equation}\vspace{0.3cm}
\hfill
$\Box$
\vspace{0.3cm}

\noindent{\bf Check of Lemma \ref{c10}}:
It follows from assumption {\bf (A.2)} that
\begin{align}
	-F_3(t,x)=&-g_{tt}(t)+g_t(t)a(x)+V(x)f_t(t)-2V(x)g(t)+V_x(x)h(t,x)+V(x)h_x(t,x) \notag\\
	\leq&-\varepsilon_2(\theta-1)(\theta-2)(1+t)^{\theta-3}+\varepsilon_2a(x)(\theta-1)(1+t)^{\theta-2} +\theta \varepsilon_1V(x)(1+t)^{\theta-1}\notag \\
	&-2\varepsilon_2V(x)(1+t)^{\theta-1}+\varepsilon_3V_x(x)x(1+t)^{\theta-1}\phi(x)+\varepsilon_3 V(x)(1+t)^{\theta-1}  \notag\\
	\leq &\varepsilon_2a(x)(1+t)^{\theta-2}\big \{ \frac{1+R+t}{a_{1}(1+t)}|(\theta-1) (\theta-2)|+|\theta-1| \big \}\notag \\
	&+V(x)(1+t)^{\theta-1}\big\{\theta \varepsilon_1-2\varepsilon_2+\varepsilon_3\big \}. \label{c33}
\end{align}
Let $\theta \varepsilon_1-2\varepsilon_2+\varepsilon_3\leq 0$, that is, $\varepsilon_3\leq 2\varepsilon_2-\theta \varepsilon_1$. At this stage, it must be $2\varepsilon_2>\theta \varepsilon_1$ to guarantee the positivity of $\varepsilon_3$. Observing \eqref{c29}, it is necessary to choose 
\begin{equation}\label{c31}
	\frac{2a_1-\theta}{2}\varepsilon_1> \varepsilon_2>\frac{\theta \varepsilon_1}{2}~~\Rightarrow~~\theta< a_1.	
\end{equation}
Note that from the assumption $0<a_1\leq 2$, and \eqref{c31} we find that $\theta < 2$. Therefore, by \eqref{c33}, Lemma \ref{c10} can be obtained when $t \gg 1$.
\noindent
\eqref{c30} and \eqref{c31} implies 
\[\theta <a_1.\]
Therefore, we can choose $\theta = a_1-\delta$ for any $\delta > 0$.
\hfill
$\Box$
\\
\vspace{0.3cm}
From the above discussion, in order to see that Lemmas \ref{c9} and \ref{c10} hold true for $\alpha=1$ and $0<a_1\leq 2$, it suffices to choose $\varepsilon_1$, 
$\varepsilon_2$, $\varepsilon_3$ and $k$ as following:\\
\[\varepsilon_1=\mu,~~~\varepsilon_2=\frac{a_1 \mu}{2},~~~ \varepsilon_3=\frac{\gamma \mu}{2},~~~ k>\max \big\{\frac{\varepsilon_3 \|a\|_\infty}{2\varepsilon_2-\theta \varepsilon_1},~\Vert a\Vert_{\infty} \big\}= \Vert a\Vert_{\infty},
\]
where $\gamma=a_1-\theta >0$ and  $\mu>0$.\\
On the other hand, for Lemma  \ref{c18}, we need to change \eqref{c15} slightly by  
\[
f(t)- |h(t,x)|
\geq  (1+t)^\theta(\varepsilon_1- \frac{\varepsilon_3}{1+t})\geq \frac{\varepsilon_1}{2} (1+t)^\theta=\frac{1}{2}f(t).\]

The rest of the proof is similar to the parts previously done. By proceeding similar arguments, we will obtain the results of the case (1) of Theorem \ref{c21-1} for $0<a_1\leq 2$ and $\alpha = 1$.

\vspace{0.5cm}
\noindent{\em Acknowledgement.}
\smallskip
This paper was written during Xiaoyan Li's stay as an overseas researcher at Hiroshima University from 12 December, 2022 to 11 December, 2023 under Ikehata's supervision as a host researcher. This work of the first author (Xiaoyan Li) was financially supported in part by Chinese Scholarship Council (Grant No. 202206160071). The work of the second author (Ryo Ikehata) was supported in part by Grant-in-Aid for Scientific Research (C) 20K03682 of JSPS.



\begin{thebibliography}{99}

\bibitem{A} L. Aloui, S. Ibrahim and M. Khenissi, Energy dcay for linear dissipative wave equation in exterior domains, J. Diff. Eq. 259 (2015), 2061-2079.
\bibitem{BR} J.-M. Bouclet and J. Royer, Local energy decay for the damped wave equation, J. Funct. Anal. 266 (2014), 4538--4615.
\bibitem{D} M. Daoulatli, Energy dcay rates for solutions of the wave equation with linear damping in exterior domain, Evol. Equ. Control Theory 5 (2016), no. 1, 37-59.
\bibitem{G} V. Georgiev, H. Kubo and K. Wakasa, Critical exponent for nonlinear damped wave equations with nonnegative potential in $3D$, J. Diff. Eq. 267 (2019), 3271--3288.
\bibitem{ikawa}M. Ikawa, Hyperbolic Partial Differential Equations and Wave Phenomena; 2000. Translations of Mathematical Monographs, American Mathematical Society.
\bibitem{I} R. Ikehata, Fast decay of solutions for linear wave equations with dissipation localized near infinity in an exterior domain, J. Diff. Eq. 188 (2003), 390-405.
\bibitem{I-1} R. Ikehata, A role of potential on $L^{2}$-estimates for some evolution equations, arXiv: 2211. 03389v1 [math AP] 7 Nov 2022.
\bibitem{II}R. Ikehata and Y. Inoue, Total energy decay for semilinear wave equations with a critical potential type of damping, Nonlinear Anal. 69 (2008), no. 4, 1396--1401.
\bibitem{IK}R. Ikehata and T. Komatsu, Fast energy decay for wave equations with variable damping coefficients in the 1-D half line. Differential Integral Equations 29 (2016), no. 5-6, 421-440.
\bibitem{IL}R. Ikehata and X. Li, Eergy decay for wave equations with a potential and a localized damping, arXiv: 2302.08114v1 [math AP], 2023.
\bibitem{IM} R. Ikehata and T. Matsuyama, $L^{2}$-behaviour of solutions to the linear heat and wave equations in exterior domains, Sci. Math. Japon. 55 (2002), 33-42.
\bibitem{IT} R. Ikehata and H. Takeda, Uniform energy decay for wave equations with unbounded damping coefficients, Funkcial. Ekvac. 63 (2020), no. 1, 133–152.
\bibitem{ITY}R. Ikehata, G. Todorova and B. Yordanov, Optimal decay rate of the energy for wave equations with critical potential, J. Math. Soc. Japan 65 (2013), no. 1, 183--236.
\bibitem{J} R. Joly and J. Royer, Energy decay and diffusion phenomenon for the asymptotically periodic damped wave equation, J. Math. Soc. Japan 70, No. 4 (2018), 1375-1418.
\bibitem{KNW} M. Kadowaki, H. Nakazawa and K. Watanabe, Exponential decay and spectral structure for wave equation with some dissipations, Tokyo J. Math. 28 (2005), no. 2, 463-470.
\bibitem{lai} N. An Lai, M. Liu, Z. Tu and C. Wang, Lifespan estimates for semilinear wave equations with space dependent damping and potential, Calc. Var. (2023), 62:44. 
\bibitem{M} A. Matsumura, On the asymptotic behavior of solutions of semi-linear wave equations, Publ. RIMS, Kyoto Univ. 12 (1976), 169--189. 
\bibitem{M-2} K. Mochizuki, Scattering theory for wave equations with dissipative terms, Publ. RIMS 12 (1976), 383-390.
\bibitem{Mo} K. Mochizuki and H. Nakazawa, Energy decay and asymptotic behavior of solutions to the wave equations with linear dissipation, Publ. RIMS 32 (1996), 401-414.
\bibitem{N}M. Nakao, Stabilization of local energy in an exterior domain for the wave equation with a localized dissipation, J. Diff. Eq. 148 (1998), no. 2, 388--406.
\bibitem{Nakao-1} M. Nakao, Energy decay for the linear and semilinear wave equations in exterior domains with some localized dissipations, Math. Z. 238 (2001), 781-797.
\bibitem{Nishihara} K. Nishihara, $Lp$-$Lq$ estimates of solutions to the damped wave equation in $3$-dimensional space and their application, Math. Z. 244 (2003), no. 3, 631--649. 
\bibitem{Nishi} H. Nishiyama, Remarks on the asymptotic behavior of the solution to damped wave equations, J. Diff. Eq. 261 (2016), no. 7, 3893--3940.
\bibitem{OZP} R. Orive, E. Zuazua and A. F. Pazoto, Asymptotic expansion for damped wave equations with periodic coefficients, Math. Models Methods Appl. Sci. 11 (2001), 1285--1310.
\bibitem{RTY} P. Radu, G. Todorova and B. Yordanov, Decay estimates for wave equations with variable coefficients, Trans. Amer. Math. Soc. 362 (2010), no. 5, 2279--2299.
\bibitem{RTY-2} P. Radu, G. Todorova and B. Yordanov, Higher order energy decay rates for damped wave equations with variable coefficients, Discrete Contin. Dyn. Syst. Ser. S 2 (2009), no. 3, 609--629.
\bibitem{RTY-3} P. Radu, G. Todorova and B. Yordanov, The generalized diffusion phenomenon and applications, SIAM J. Math. Anal. 48 (2016), no. 1, 174--203.
\bibitem{S}M. Sobajima, On global existence for semilinear wave equations with space-dependent critical damping, J. Math. Soc. Japan 75, No.2 (2023), 603--627.
\bibitem{S-1}M. Sobajima, Asymptotic behavior for wave equations with space-dependent damping in a weighted energy class, Communications on Pure Appl. Anal. (to appear).
\bibitem{SW} M. Sobajima and Y. Wakasugi, Diffusion phenomena for the wave equation with space-dependent damping in an exterior domain, J. Diff. Eq. 261 (2016), no. 10, 5690--5718.
\bibitem{SW-1} M. Sobajima and Y. Wakasugi, Remarks on an elliptic problem arising in weighted energy estimates for wave equations with space-dependent damping term in an exterior domain, AIMS Math. 2:1-15, 2017.
\bibitem{SW-2} M. Sobajima and Y. Wakasugi, Diffusion phenomena for the wave equation with space-dependent damping term growing at infinity, Adv. Differential Equations 23 (2018), no. 7-8, 581--614.
\bibitem{TY} G. Todorova and B. Yordanov, Weighted $L2$-estimates of dissipative wave equations with variable coefficients, J. Diff. Eq. 246 (2009), no. 12, 4497--4518.
\bibitem{W}Y. Wakasugi, On diffusion phenomena for the linear wave equation with space-dependent damping, J. Hyperblic Differ. Equ. 11(2014), 795-819.
\bibitem{W-2}Y. Wakasugi, Decay property of solutions to the wave equation with space-dependent damping, absorbing nonlinearity, and polynomially decaying data, Math. Meth. Appl. Sci. 2023;1-41. doi:10.1002/mma.8957
\bibitem{Z} Z. Zhang, Fast decay of solutions for wave equations with localized dissipation on noncompact Riemannian manifold, Nonlinear Anal. Real World Appl. 2 (2016), 246-260.
\bibitem{Zua}E. Zuazua, Exponential decay for the semilinear wave equation with localized damping in unbounded domains, J. Math. Pures et Appl. 70 (1992), 513-529.
 

\end{thebibliography}
\end{document}